\documentclass[12pt,reqno]{article}

\usepackage[usenames]{color}
\usepackage[colorlinks=true,
linkcolor=webgreen, filecolor=webbrown,
citecolor=webgreen]{hyperref}

\definecolor{webgreen}{rgb}{0,.5,0}
\definecolor{webbrown}{rgb}{.6,0,0}

\usepackage{amssymb}
\usepackage{graphicx}
\usepackage{amscd}
\usepackage{lscape}
\usepackage{tikz}
\usepackage{tikz-cd}
\usepackage{pgfplots}

\usetikzlibrary{matrix}
\usetikzlibrary{fit,shapes}
\usetikzlibrary{positioning, calc}
\tikzset{circle node/.style = {circle,inner sep=1pt,draw, fill=white},
        X node/.style = {fill=white, inner sep=1pt},
        dot node/.style = {circle, draw, inner sep=5pt}
        }
\usepackage{tkz-fct}

\usepackage{amsthm}
\newtheorem{theorem}{Theorem}
\newtheorem{lemma}[theorem]{Lemma}
\newtheorem{proposition}[theorem]{Proposition}
\newtheorem{corollary}[theorem]{Corollary}

\theoremstyle{definition}

\newtheorem{example}[theorem]{Example}

\usepackage{float}

\usepackage{graphics,amsmath}
\usepackage{amsfonts}
\usepackage{latexsym}
\usepackage{epsf}

\newcommand{\seqnum}[1]{\href{http://oeis.org/#1}{\underline{#1}}}

\begin{document}

\begin{center}
\vskip 1cm{\LARGE\bf A note on number triangles that are almost their own production matrix} \vskip 1cm \large
Paul Barry\\
School of Science\\
Waterford Institute of Technology\\
Ireland\\
\href{mailto:pbarry@wit.ie}{\tt pbarry@wit.ie}
\end{center}
\vskip .2 in

\begin{abstract} We characterize a family of number triangles whose production matrices are closely related to the original number triangle. We study a number of such triangles that are of combinatorial significance. For a specific subfamily, these triangles relate to sequences that have interesting convolution recurrences and continued fraction generating functions.
\end{abstract}

\section{Preliminaries} In this note, we shall be concerned with infinite lower-triangular matrices with integer entries. Thus we will be working with matrices $(T_{n,k})$ where $T_{n,k}=0$ if $k>n$. We shall index these matrices starting with $T_{0,0}$ in the upper left corner. The Riordan group of ordinary Riordan arrays \cite{Book, Survey, SGWW} is composed of such matrices, where an ordinary Riordan array is defined by two power series,
$$g(x)=1+g_1 x+ g_2 x^2+\cdots, \quad f(x)=x+f_2x^2+f_3x^3+\cdots,$$ where 
$$T_{n,k}=[x^n] g(x)f(x)^k.$$ 
Here, $[x^n]$ is the operator that extracts the coefficient of the term in $x^n$. 
Similarly, an exponential Riordan array $[g(x),f(x)]$ \cite{Book, DeutschShap} is defined by two power series 
$$g(x)=1+g_1 \frac{x}{1!}+ g_2 \frac{x^2}{2!}+\cdots, \quad f(x)=\frac{x}{1!}+f_2\frac{x^2}{2!}+f_3\frac{x^3}{3!}+\cdots,$$ where 
the general $(n,k)$-th element of $[g(x),f(x)]$ is given by 
$$T_{n,k}=\frac{n!}{k!}[x^n]g(x)f(x)^k.$$ 

In dealing with the infinite matrix $(T_{n,k})$, we display a suitable truncation in the text that follows. 

Many of the sequences that we shall encounter have ordinary generating functions that can be expressed as Stieltjes or Jacobi continued fractions \cite{CFT, Wall}. For these we use the following notation. 
The Stieltjes continued fraction
$$\cfrac{1}{1-
\cfrac{ax}{1-
\cfrac{\alpha x}{1-
\cfrac{bx}{1-
\cfrac{\beta x}{1-\cdots}}}}}$$ will be denoted by $\mathcal{S}(a,\alpha,b,\beta,\ldots)$ or where appropriate as $\mathcal{S}(a,b,\ldots;\alpha,\beta,\ldots)$.

The Jacobi continued fraction 
$$\cfrac{1}{1- ax-
\cfrac{bx^2}{1-cx-
\cfrac{dx^2}{1-\cdots}}}$$ will be denoted by $\mathcal{J}(a,c,\ldots;b,d,\ldots)$. 

Sequences and triangles, where known, will be referenced by their $Annnnnn$ number in the On-Line Encyclopedia of Integer Sequences \cite{SL1, SL2}.

\section{Introduction} For an invertible lower-triangular matrix $M$, we define its production matrix \cite{ProdMat} to be the matrix $$P_M= M^{-1} \overline{M}= M^{-1}\cdot U\cdot M,$$ where $U$ denotes the infinite shift matrix that begins
$$\left(
\begin{array}{ccccccc}
 0 & 1 & 0 & 0 & 0 & 0 & 0 \\
 0 & 0 & 1 & 0 & 0 & 0 & 0 \\
 0 & 0 & 0 & 1 & 0 & 0 & 0 \\
 0 & 0 & 0 & 0 & 1 & 0 & 0 \\
 0 & 0 & 0 & 0 & 0 & 1 & 0 \\
 0 & 0 & 0 & 0 & 0 & 0 & 1 \\
 0 & 0 & 0 & 0 & 0 & 0 & 0 \\
\end{array}
\right).$$
Thus $\overline{M}$ is the matrix $M$ with its top row removed.

By construction, the matrix $P_M$ will have a non-zero superdiagonal, so it will differ from $M$ in this respect (note that $P_M$ is a Hessenberg matrix). For the purposes of this note, we shall say that the matrix $M$ is \emph{almost its own production matrix} $P_M$ if, apart from the extra super-diagonal in $P_M$, the two matrices only differ otherwise in their first column.

In order to define a class of number triangles that are almost their own production matrix, we introduce some notation. Given a sequence $b_0, b_1, b_2, \ldots$ with generating function $g(x)$, we define the corresponding \emph{sequence array} to be the lower triangular matrix with $(n,k)$-term equal to $b_{n-k}$, for $k \le n$, and $0$ otherwise. In the language of ordinary Riordan arrays, this is the array $(g(x),x)$, a member of the Appell subgroup.

We define $V$ to be the infinite matrix with $-n$ on the sub-diagonal, and $0$ otherwise. Thus $V$ begins
$$\left(
\begin{array}{ccccccc}
 0 & 0 & 0 & 0 & 0 & 0 & 0 \\
 0 & 0 & 0 & 0 & 0 & 0 & 0 \\
 0 & -1 & 0 & 0 & 0 & 0 & 0 \\
 0 & 0 & -2 & 0 & 0 & 0 & 0 \\
 0 & 0 & 0 & -3 & 0 & 0 & 0 \\
 0 & 0 & 0 & 0 & -4 & 0 & 0 \\
 0 & 0 & 0 & 0 & 0 & -5 & 0 \\
\end{array}
\right).$$

We have the following proposition.
\begin{proposition}
Let $a_n$ be a sequence with generating function $f(x)=a_0+a_1x+a_2 x^2+\cdots$. Then the lower-triangular matrix 
$$M=( (1-xf(x),x)+V)^{-1}$$ is almost its own production matrix. In fact, we have that
$$P_M-M$$ takes the form of the matrix that begins
$$\left(
\begin{array}{cccccc}
 a_0-1 & 1 & 0 & 0 & 0 & 0 \\
 a_1 & 0 & 1 & 0 & 0 & 0 \\
 a_2 & 0 & 0 & 1 & 0 & 0 \\
 a_3 & 0 & 0 & 0 & 1 & 0 \\
 a_4 & 0 & 0 & 0 & 0 & 1 \\
 a_5 & 0 & 0 & 0 & 0 & 0 \\
\end{array}
\right).$$
\end{proposition}

\begin{example} We let $a_n$ be the sequence of indecomposable permutations \seqnum{A003319} that begins
$$ 1, 1, 3, 13, 71, 461, 3447, 29093, 273343, 2829325, 31998903,\ldots.$$ This sequence is the INVERT$(1)$ transform of the shifted factorials $(n+1)!$. The (ordinary) generating function $f(x)$ of this sequence is given by the Stieltjes continued fraction
$$\cfrac{1}{1-
\cfrac{x}{1-
\cfrac{2x}{1-
\cfrac{2x}{1-
\cfrac{3x}{1-
\cfrac{3x}{1-\ldots}}}}}},$$ with coefficients $1,2,2,3,3,4,4,\ldots$.

The matrix $M=( (1-xf(x),x)+V)^{-1}$ then begins
$$\left(
\begin{array}{ccccccc}
 1 & 0 & 0 & 0 & 0 & 0 & 0 \\
 -1 & 1 & 0 & 0 & 0 & 0 & 0 \\
 -1 & -2 & 1 & 0 & 0 & 0 & 0 \\
 -3 & -1 & -3 & 1 & 0 & 0 & 0 \\
 -13 & -3 & -1 & -4 & 1 & 0 & 0 \\
 -71 & -13 & -3 & -1 & -5 & 1 & 0 \\
 -461 & -71 & -13 & -3 & -1 & -6 & 1 \\
\end{array}
\right)^{-1}=\left(
\begin{array}{ccccccc}
 1 & 0 & 0 & 0 & 0 & 0 & 0 \\
 1 & 1 & 0 & 0 & 0 & 0 & 0 \\
 3 & 2 & 1 & 0 & 0 & 0 & 0 \\
 13 & 7 & 3 & 1 & 0 & 0 & 0 \\
 71 & 33 & 13 & 4 & 1 & 0 & 0 \\
 461 & 191 & 71 & 21 & 5 & 1 & 0 \\
 3447 & 1297 & 461 & 133 & 31 & 6 & 1 \\
\end{array}
\right).$$
This is \seqnum{A104980} [Paul D. Hanna, \cite{SL1}]. The production matrix $P_M$ of $M$ then begins
$$\left(
\begin{array}{ccccccc}
 1 & 1 & 0 & 0 & 0 & 0 & 0 \\
 2 & 1 & 1 & 0 & 0 & 0 & 0 \\
 6 & 2 & 1 & 1 & 0 & 0 & 0 \\
 26 & 7 & 3 & 1 & 1 & 0 & 0 \\
 142 & 33 & 13 & 4 & 1 & 1 & 0 \\
 922 & 191 & 71 & 21 & 5 & 1 & 1 \\
 6894 & 1297 & 461 & 133 & 31 & 6 & 1 \\
\end{array}
\right).$$

We have that $P_M-M$ begins
$$\left(
\begin{array}{ccccccc}
 1 & 1 & 0 & 0 & 0 & 0 & 0 \\
 2 & 1 & 1 & 0 & 0 & 0 & 0 \\
 6 & 2 & 1 & 1 & 0 & 0 & 0 \\
 26 & 7 & 3 & 1 & 1 & 0 & 0 \\
 142 & 33 & 13 & 4 & 1 & 1 & 0 \\
 922 & 191 & 71 & 21 & 5 & 1 & 1 \\
 6894 & 1297 & 461 & 133 & 31 & 6 & 1 \\
\end{array}
\right)-\left(
\begin{array}{ccccccc}
 1 & 0 & 0 & 0 & 0 & 0 & 0 \\
 1 & 1 & 0 & 0 & 0 & 0 & 0 \\
 3 & 2 & 1 & 0 & 0 & 0 & 0 \\
 13 & 7 & 3 & 1 & 0 & 0 & 0 \\
 71 & 33 & 13 & 4 & 1 & 0 & 0 \\
 461 & 191 & 71 & 21 & 5 & 1 & 0 \\
 3447 & 1297 & 461 & 133 & 31 & 6 & 1 \\
\end{array}
\right)$$
$$=\left(
\begin{array}{ccccccc}
 0 & 1 & 0 & 0 & 0 & 0 & 0 \\
 1 & 0 & 1 & 0 & 0 & 0 & 0 \\
 3 & 0 & 0 & 1 & 0 & 0 & 0 \\
 13 & 0 & 0 & 0 & 1 & 0 & 0 \\
 71 & 0 & 0 & 0 & 0 & 1 & 0 \\
 461 & 0 & 0 & 0 & 0 & 0 & 1 \\
 3447 & 0 & 0 & 0 & 0 & 0 & 0 \\
\end{array}
\right).$$ 
Rearranging, we obtain that the production matrix $P_M$ can then be expressed as beginning
$$\left(
\begin{array}{ccccccc}
 1 & 1 & 0 & 0 & 0 & 0 & 0 \\
 2 & 1 & 1 & 0 & 0 & 0 & 0 \\
 6 & 2 & 1 & 1 & 0 & 0 & 0 \\
 26 & 7 & 3 & 1 & 1 & 0 & 0 \\
 142 & 33 & 13 & 4 & 1 & 1 & 0 \\
 922 & 191 & 71 & 21 & 5 & 1 & 1 \\
 6894 & 1297 & 461 & 133 & 31 & 6 & 1 \\
\end{array}
\right)=$$ 
$$\left(
\begin{array}{ccccccc}
 1 & 0 & 0 & 0 & 0 & 0 & 0 \\
 1 & 1 & 0 & 0 & 0 & 0 & 0 \\
 3 & 2 & 1 & 0 & 0 & 0 & 0 \\
 13 & 7 & 3 & 1 & 0 & 0 & 0 \\
 71 & 33 & 13 & 4 & 1 & 0 & 0 \\
 461 & 191 & 71 & 21 & 5 & 1 & 0 \\
 3447 & 1297 & 461 & 133 & 31 & 6 & 1 \\
\end{array}
\right)+\left(
\begin{array}{ccccccc}
 0 & 1 & 0 & 0 & 0 & 0 & 0 \\
 1 & 0 & 1 & 0 & 0 & 0 & 0 \\
 3 & 0 & 0 & 1 & 0 & 0 & 0 \\
 13 & 0 & 0 & 0 & 1 & 0 & 0 \\
 71 & 0 & 0 & 0 & 0 & 1 & 0 \\
 461 & 0 & 0 & 0 & 0 & 0 & 1 \\
 3447 & 0 & 0 & 0 & 0 & 0 & 0 \\
\end{array}
\right).$$ 
\end{example}

\section{Proof of the proposition}
The proof of the proposition hinges on two lemmas. 
\begin{lemma} The production matrix of an element $(g(x),x)^{-1}$ of the Appell subgroup of the Riordan group, where $g(x)=1+b_1x+b_2x+ \ldots$,  begins
$$\left(
\begin{array}{ccccccc}
 -b_1 & 1 & 0 & 0 & 0 & 0 & 0 \\
 -b_2 & 0 & 1 & 0 & 0 & 0 & 0 \\
 -b_3 & 0 & 0 & 1 & 0 & 0 & 0 \\
 -b_4 & 0 & 0 & 0 & 1 & 0 & 0 \\
 -b_5 & 0 & 0 & 0 & 0 & 1 & 0 \\
 -b_6 & 0 & 0 & 0 & 0 & 0 & 1 \\
 -b_7 & 0 & 0 & 0 & 0 & 0 & 0 \\
\end{array}
\right).$$ 
\end{lemma}
\begin{corollary} The production matrix of the element $(1-xf(x),x)^{-1}$ of the Appell subgroup, where 
$f(x)=1+a_1 x+ a_2 x + \cdots$, begins
$$\left(\begin{array}{ccccccc}
 a_0 & 1 & 0 & 0 & 0 & 0 & 0 \\
 a_1 & 0 & 1 & 0 & 0 & 0 & 0 \\
 a_2 & 0 & 0 & 1 & 0 & 0 & 0 \\
 a_3 & 0 & 0 & 0 & 1 & 0 & 0 \\
 a_4 & 0 & 0 & 0 & 0 & 1 & 0 \\
 a_5 & 0 & 0 & 0 & 0 & 0 & 1 \\
 a_6 & 0 & 0 & 0 & 0 & 0 & 0 \\
\end{array}
\right).$$ 
\end{corollary}
For the next lemma, we let $W$ be the infinite matrix with $1$ in its top left position, and $0$ elsewhere. The identity matrix will be denoted by $I$. 
\begin{lemma} We have the following identity of matrices.
$$V \cdot U - U \cdot V = I-W.$$
\end{lemma}

We can now prove the proposition of the previous section. We are required to prove that 
$$P_M-M=P_{A^{-1}}-W,$$ where $A=(1-xf(x),x)$. We shall let $A_0$ denote the matrix whose first column is that of $A$, and all of whose other entries are $0$. Thus we have 
$$P_{A^{-1}}=A_0+U.$$ 
The desired equality 
$$P_M-M=P_{A^{-1}}-W$$ can now be written 
$$(A+V)U(A+V)^{-1}-(A+V)^{-1}=A_0+U-W.$$ 
Multiply both sides by $A+V$ on the right, to obtain an equivalent statement that must be proved. We get
$$(A+V)U-I=A_0(A+V)+U(A+V)-W(A+V),$$ which is equivalent to 
$$(A+V)U-I=A_0+U(A+V)-W.$$ Equivalently, we require that 
$$AU+VU-I=A_0+UA+UV-W.$$ But this is so since $V  U - U  V = I-W$. 
The result is thus proven.

\begin{example}
We let $f(x)=\frac{1}{1-x}$, so that the sequence $a_n$ begins 
$$1,1,1,1,1,\ldots.$$ 
We find that 
$$M=\left(
\begin{array}{ccccccc}
 1 & 0 & 0 & 0 & 0 & 0 & 0 \\
 1 & 1 & 0 & 0 & 0 & 0 & 0 \\
 3 & 2 & 1 & 0 & 0 & 0 & 0 \\
 11 & 7 & 3 & 1 & 0 & 0 & 0 \\
 49 & 31 & 13 & 4 & 1 & 0 & 0 \\
 261 & 165 & 69 & 21 & 5 & 1 & 0 \\
 1631 & 1031 & 431 & 131 & 31 & 6 & 1 \\
\end{array}
\right),$$ and 
that 
$$P_M=\left(
\begin{array}{ccccccc}
 1 & 1 & 0 & 0 & 0 & 0 & 0 \\
 2 & 1 & 1 & 0 & 0 & 0 & 0 \\
 4 & 2 & 1 & 1 & 0 & 0 & 0 \\
 12 & 7 & 3 & 1 & 1 & 0 & 0 \\
 50 & 31 & 13 & 4 & 1 & 1 & 0 \\
 262 & 165 & 69 & 21 & 5 & 1 & 1 \\
 1632 & 1031 & 431 & 131 & 31 & 6 & 1 \\
\end{array}
\right).$$ 
Thus we have 
$$P_M-M=\left(
\begin{array}{ccccccc}
 0 & 1 & 0 & 0 & 0 & 0 & 0 \\
 1 & 0 & 1 & 0 & 0 & 0 & 0 \\
 1 & 0 & 0 & 1 & 0 & 0 & 0 \\
 1 & 0 & 0 & 0 & 1 & 0 & 0 \\
 1 & 0 & 0 & 0 & 0 & 1 & 0 \\
 1 & 0 & 0 & 0 & 0 & 0 & 1 \\
 1 & 0 & 0 & 0 & 0 & 0 & 0 \\
\end{array}
\right).$$ 
The self-building nature of these triangles is illustrated below: For any element of a non-initial column, the elements above it in the same column, prepended by $1$, are used as multipliers on the elements in the row above it, starting in the column to its left.
\begin{center}
\begin{tikzpicture}
\matrix[matrix of math nodes,left delimiter = (,right delimiter = ),row sep=2pt,column sep = 2pt] (o){
 1 & 0 & 0 & 0 & 0 & 0 & 0 \\
 1 & 1 & 0 & 0 & 0 & 0 & 0 \\
 {}_{}^1 3 & {}_{}^2 2 &{}_{}^4 1 & 0 & 0 & 0 & 0 \\
 11 & 7 & 3 & 1 & 0 & 0 & 0 \\
  {}_{}^1 49 &  {}_{}^1 31 &  {}_{}^2 13 &{}_{}^7 4 &{}_{}^{31} 1 & 0 & 0 \\
 261 & 165 & 69 & {}_{}^1 21 & {}_{}^1 5 & {}_{}^{5} 1 & 0 \\
 1631 & 1031 & 431 & 131 & 31 & 6 & 1 \\
};
\draw (o-3-1.north west) rectangle (o-3-3.south east);
\draw (o-4-1.north west) rectangle (o-4-1.south east);
\draw (o-5-1.north west) rectangle (o-5-5.south east);
\draw (o-6-2.north west) rectangle (o-6-2.south east);
\draw (o-6-4.north west) rectangle (o-6-6.south east);
\draw (o-7-5.north west) rectangle (o-7-5.south east);
\end{tikzpicture}
\end{center}
For the initial column, we use the elements in the first column of $P_M$  in a similar fashion, only this time starting in the row immediately above. (Note that for consistency, we could have prepended a $1$ to this sequence, and used it on an all-zero column to the left of the initial column of the matrix).

The construction is thus as follows. We start with a sequence 
$$a_0, a_1, a_2, \ldots.$$ 
We set $T_{0,0}=a_0$, and $T_{n,k}=0$ if $k>n$. Then we have 
$$T_{n,0}=(a_0+m_0-1) T_{n-1,0}+ \sum_{j=1}^{n-1}(a_j+m_j) T_{n-1,j},$$ 
and for $k>0$, 
$$T_{n,k}=T_{n-1,k-1}+\sum_{j=0}^{n-k-1} T_{k+j,k}T_{n-1,k+j}.$$ 

Here, the sequence $m_0, m_1, \ldots$ is the first column of $M$, and the sequence $a_n+m_n-0^n$ is the first column of the production matrix $P_M$. 

Note that the sequence $1,3,11,49,\ldots$ \seqnum{A001339} is the binomial transform of $(n+1)!$, with exponential generating function $\frac{e^x}{(1-x)^2}$. Its ordinary generating function is given by the Jacobi continued fraction $\mathcal{J}(3,5,7,9,\ldots;2,6,12,20,\ldots)$. 
\end{example}

\begin{example} We set $f(x)=\frac{1}{1-2x}$, so that the sequence $a_0, a_1, a_2,\ldots$ is the sequence $a_n=2^n$. We find that 
$$M=\left(
\begin{array}{ccccccc}
 1 & 0 & 0 & 0 & 0 & 0 & 0 \\
 1 & 1 & 0 & 0 & 0 & 0 & 0 \\
 4 & 2 & 1 & 0 & 0 & 0 & 0 \\
 18 & 8 & 3 & 1 & 0 & 0 & 0 \\
 92 & 40 & 14 & 4 & 1 & 0 & 0 \\
 536 & 232 & 80 & 22 & 5 & 1 & 0 \\
 3552 & 1536 & 528 & 144 & 32 & 6 & 1 \\
\end{array}
\right),$$ and 
$$P_M=\left(
\begin{array}{ccccccc}
 1 & 1 & 0 & 0 & 0 & 0 & 0 \\
 3 & 1 & 1 & 0 & 0 & 0 & 0 \\
 8 & 2 & 1 & 1 & 0 & 0 & 0 \\
 26 & 8 & 3 & 1 & 1 & 0 & 0 \\
 108 & 40 & 14 & 4 & 1 & 1 & 0 \\
 568 & 232 & 80 & 22 & 5 & 1 & 1 \\
 3616 & 1536 & 528 & 144 & 32 & 6 & 1 \\
\end{array}
\right).$$ 
Thus we have 
$$P_M-P=\left(
\begin{array}{ccccccc}
 0 & 1 & 0 & 0 & 0 & 0 & 0 \\
 2 & 0 & 1 & 0 & 0 & 0 & 0 \\
 4 & 0 & 0 & 1 & 0 & 0 & 0 \\
 8 & 0 & 0 & 0 & 1 & 0 & 0 \\
 16 & 0 & 0 & 0 & 0 & 1 & 0 \\
 32 & 0 & 0 & 0 & 0 & 0 & 1 \\
 64 & 0 & 0 & 0 & 0 & 0 & 0 \\
\end{array}
\right).$$ 
The sequence $1, 4, 18, 92, 536, 3552, \ldots$ \seqnum{A081923} has exponential generating function $\frac{e^{2x}}{(1-x)^2}$, and is the binomial generating function of the sequence in the previous example. Its ordinary generating function is thus given by $\mathcal{J}(4,6,8,10,\ldots;2,6,12,20,\ldots)$. 
\end{example}
Following from these two examples, we have the general result that if $f(x)=\frac{1}{1-rx}$, and hence $a_n=r^n$, then the first column of $M$ is given by 
$$m_n=\sum_{k=0}^n\binom{n-1}{k}(n-k)!r^k=\sum_{k=0}^n \binom{n-1}{n-k}r^{n-k}k!,$$ with coefficient array that begins
$$\left(
\begin{array}{ccccccc}
 1 & 0 & 0 & 0 & 0 & 0 & 0 \\
 1 & 0 & 0 & 0 & 0 & 0 & 0 \\
 2 & 1 & 0 & 0 & 0 & 0 & 0 \\
 6 & 4 & 1 & 0 & 0 & 0 & 0 \\
 24 & 18 & 6 & 1 & 0 & 0 & 0 \\
 120 & 96 & 36 & 8 & 1 & 0 & 0 \\
 720 & 600 & 240 & 60 & 10 & 1 & 0 \\
\end{array}
\right).$$ 
This is essentially \seqnum{A132159}. 
The sequence $m_{n+1}$ then has exponential generating sequence $\frac{e^{rx}}{(1-x)^2}$ and ordinary generating function $\mathcal{J}(r+2,r+4,r+6,10,\ldots;2,6,12,20,\ldots)$. It is thus the $r$-th binomial transform of $(n+1)!$, with $m_{n+1}=\sum_{k=0}^n \binom{n}{k}r^{n-k}(k+1)!$

\section{A family of Hanna triangles}
Paul D. Hanna has defined a one-parameter family of triangles that are almost their own production matrix, one of which we have already seen. For a parameter $r \in \mathbb{N}$, these triangles $H(r)$ can be defined as
follows: if $n<0$ or $k<0$, then $H_{n,k}=0$; if $n=k$ then $H_{n,k}=1$; if $n=k+1$ then $H_{n,k}=n$; otherwise we have
$$H_{n,k}=k H_{n,k+1}+\sum_{j=0}^{n-k-1} H_{j+r,r} H_{n,j+k+1}.$$
For $r=1,2,3,4$, we get the following triangles.
$$\left(
\begin{array}{ccccccc}
 1 & 0 & 0 & 0 & 0 & 0 & 0 \\
 1 & 1 & 0 & 0 & 0 & 0 & 0 \\
 3 & 2 & 1 & 0 & 0 & 0 & 0 \\
 13 & 7 & 3 & 1 & 0 & 0 & 0 \\
 71 & 33 & 13 & 4 & 1 & 0 & 0 \\
 461 & 191 & 71 & 21 & 5 & 1 & 0 \\
 3447 & 1297 & 461 & 133 & 31 & 6 & 1 \\
\end{array}
\right),\left(
\begin{array}{ccccccc}
 1 & 0 & 0 & 0 & 0 & 0 & 0 \\
 1 & 1 & 0 & 0 & 0 & 0 & 0 \\
 4 & 2 & 1 & 0 & 0 & 0 & 0 \\
 22 & 8 & 3 & 1 & 0 & 0 & 0 \\
 148 & 44 & 14 & 4 & 1 & 0 & 0 \\
 1156 & 296 & 84 & 22 & 5 & 1 & 0 \\
 10192 & 2312 & 600 & 148 & 32 & 6 & 1 \\
\end{array}
\right),$$
$$\left(
\begin{array}{ccccccc}
 1 & 0 & 0 & 0 & 0 & 0 & 0 \\
 1 & 1 & 0 & 0 & 0 & 0 & 0 \\
 5 & 2 & 1 & 0 & 0 & 0 & 0 \\
 33 & 9 & 3 & 1 & 0 & 0 & 0 \\
 261 & 57 & 15 & 4 & 1 & 0 & 0 \\
 2361 & 441 & 99 & 23 & 5 & 1 & 0 \\
 23805 & 3933 & 783 & 165 & 33 & 6 & 1 \\
\end{array}
\right),\left(
\begin{array}{ccccccc}
 1 & 0 & 0 & 0 & 0 & 0 & 0 \\
 1 & 1 & 0 & 0 & 0 & 0 & 0 \\
 6 & 2 & 1 & 0 & 0 & 0 & 0 \\
 46 & 10 & 3 & 1 & 0 & 0 & 0 \\
 416 & 72 & 16 & 4 & 1 & 0 & 0 \\
 4256 & 632 & 116 & 24 & 5 & 1 & 0 \\
 48096 & 6352 & 1016 & 184 & 34 & 6 & 1 \\
\end{array}
\right).$$ 
These are \seqnum{A104980}, \seqnum{A111536}, \seqnum{111544}, \seqnum{A111553} respectively. The corresponding first column elements are \seqnum{A003319}, \seqnum{A111529}, \seqnum{A111530} and \seqnum{A111531}, respectively.

The defining sequences $a_0, a_1, a_2,\ldots$ of these triangles can be characterised as follows. For a given $r \in \mathbb{N}$, the sequence $a_0, a_1, a_2,\ldots$ is the sequence with generating function given by the Stieltjes continued fraction
$$\mathcal{S}(r,r+1,r+2,\ldots;2,3,4,\ldots).$$
Thus when $r=0$, we have $f(x)=1$, and we obtain the triangle \seqnum{A094587}
$$\left(
\begin{array}{ccccccc}
 1 & 0 & 0 & 0 & 0 & 0 & 0 \\
 1 & 1 & 0 & 0 & 0 & 0 & 0 \\
 2 & 2 & 1 & 0 & 0 & 0 & 0 \\
 6 & 6 & 3 & 1 & 0 & 0 & 0 \\
 24 & 24 & 12 & 4 & 1 & 0 & 0 \\
 120 & 120 & 60 & 20 & 5 & 1 & 0 \\
 720 & 720 & 360 & 120 & 30 & 6 & 1 \\
\end{array}
\right)=\left(
\begin{array}{ccccccc}
 1 & 0 & 0 & 0 & 0 & 0 & 0 \\
 -1 & 1 & 0 & 0 & 0 & 0 & 0 \\
 0 & -2 & 1 & 0 & 0 & 0 & 0 \\
 0 & 0 & -3 & 1 & 0 & 0 & 0 \\
 0 & 0 & 0 & -4 & 1 & 0 & 0 \\
 0 & 0 & 0 & 0 & -5 & 1 & 0 \\
 0 & 0 & 0 & 0 & 0 & -6 & 1 \\
\end{array}
\right)^{-1},$$ with production matrix
$$\left(
\begin{array}{ccccccc}
 1 & 1 & 0 & 0 & 0 & 0 & 0 \\
 1 & 1 & 1 & 0 & 0 & 0 & 0 \\
 2 & 2 & 1 & 1 & 0 & 0 & 0 \\
 6 & 6 & 3 & 1 & 1 & 0 & 0 \\
 24 & 24 & 12 & 4 & 1 & 1 & 0 \\
 120 & 120 & 60 & 20 & 5 & 1 & 1 \\
 720 & 720 & 360 & 120 & 30 & 6 & 1 \\
\end{array}
\right).$$
This is the exponential Riordan array $\left[\frac{1}{1-x},x\right]$. In this special case, we have 
$$P_M-M=U.$$ 
Returning to the general case, the sequences $m_n$ (the initial columns of the triangles $H(r)$ in question)  have their generating function given by
$$\mathcal{S}(1,2,3,\ldots; r+1,r+2,r+2,\ldots).$$ 
These sequences are solutions to the convolution recurrence 
$$m_n=(n-r)m_{n-1}+r \sum_{i=0}^{n-1} m_i m_{n-i-1},$$ with $m_0=1$, $m_1=1$. 
For general $r$, this sequence begins 
$$1, 1, r + 2, r^2 + 6r + 6, r^3 + 12r^2 + 34r + 24, r^4 + 20r^3 + 110r^2 + 210r + 120,\ldots,$$ with coefficient array that begins
$$\left(
\begin{array}{cccccccc}
 1 & 0 & 0 & 0 & 0 & 0 & 0 & 0 \\
 1 & 0 & 0 & 0 & 0 & 0 & 0 & 0 \\
 2 & 1 & 0 & 0 & 0 & 0 & 0 & 0 \\
 6 & 6 & 1 & 0 & 0 & 0 & 0 & 0 \\
 24 & 34 & 12 & 1 & 0 & 0 & 0 & 0 \\
 120 & 210 & 110 & 20 & 1 & 0 & 0 & 0 \\
 720 & 1452 & 974 & 270 & 30 & 1 & 0 & 0 \\
 5040 & 11256 & 8946 & 3248 & 560 & 42 & 1 & 0 \\
\end{array}
\right).$$ 
In the Del\'eham notation, this is a variant of the triangle
$$[0, 2, 1, 3, 2, 4, 3, 5, \ldots] \,\Delta \,[1, 0, 1, 0, 1, 0,1,0,\ldots]$$ 
which is \seqnum{A111184}.
For $r=1$, we illustrate below the ``almost its own production matrix'' property of $H(1)$. 
\begin{center}
\begin{tikzpicture}
\matrix[matrix of math nodes,left delimiter = (,right delimiter = ),row sep=2pt,column sep = 2pt] (o){
 1 & 0 & 0 & 0 & 0 & 0 & 0 \\
 1 & 1 & 0 & 0 & 0 & 0 & 0 \\
 {}_{}^1 3 & {}_{}^2 2 &{}_{}^6 1 & 0 & 0 & 0 & 0 \\
 13 & 7 & 3 & 1 & 0 & 0 & 0 \\
  {}_{}^1 71 &  {}_{}^1 33 &  {}_{}^2 13 &{}_{}^7 4 &{}_{}^{33} 1 & 0 & 0 \\
 461 & 191 & 71 & {}_{}^1 21 & {}_{}^1 5 & {}_{}^{5} 1 & 0 \\
 3447 & 1297 & 461 & 133 & 31 & 6 & 1 \\
};
\draw (o-3-1.north west) rectangle (o-3-3.south east);
\draw (o-4-1.north west) rectangle (o-4-1.south east);
\draw (o-5-1.north west) rectangle (o-5-5.south east);
\draw (o-6-2.north west) rectangle (o-6-2.south east);
\draw (o-6-4.north west) rectangle (o-6-6.south east);
\draw (o-7-5.north west) rectangle (o-7-5.south east);
\end{tikzpicture}
\end{center}

The sequences $a_n$ and $m_n$ of this section actually belong to a single family of sequences that are specified in general by three parameters (two of these parameters are equal to $1$ in our case). These have been studied by Martin and Kearney \cite{Martin1, Martin2}.
If we change the indexing of these sequences from $0$ to $1$, then $a_n$ corresponds to 
$$S(1,r-3,1)$$ 
using the notation of \cite{Martin1}, while $m_n$ corresponds to 
$$S(1,-(r+1), 1).$$ 
When $r=1$, we obtain $a_n=m_n$, the sequence of indecomposable permutations. 

The sequence $S(\alpha, \beta, \gamma)=(u_n)_{n=1}^{\infty}$ satisfies the recurrence 
$$u_n=(\alpha n + \beta) u_{n-1}+\gamma \sum_{j=1}^{n-1} u_j u_{n-j},$$ with 
$u_1=1$. 

We posit that these moment sequences have ordinary generating function 
$$\mathcal{S}(2\alpha+\beta+\gamma, 3\alpha+\beta+\gamma, 4\alpha+\beta+\gamma,\ldots;\alpha+\gamma, 2\alpha+\gamma, 3\alpha+\gamma,\ldots),$$ or equivalently,
$$\mathcal{J}(2\alpha+\beta+\gamma,4\alpha+\beta+2\gamma,6\alpha+\beta+2\gamma,\ldots;
(\alpha+\gamma)(2\alpha+\beta+\gamma), (2\alpha+\beta)(3\alpha+\beta+\gamma),(3\alpha+\beta)(4\alpha+\beta+\gamma),\ldots).$$

\bigskip
\hrule
\bigskip
\noindent 2010 {\it Mathematics Subject Classification}: Primary 15B36;
 Secondary 05A15, 11C20, 11B37, 11B83, 15B36.
\noindent \emph{Keywords:}  Number triangle, production matrix, Riordan array, convolution recurrence, moment sequence.

\bigskip
\hrule
\bigskip
\noindent (Concerned with sequences

\seqnum{A003319},
\seqnum{A001339},
\seqnum{A081923},
\seqnum{A094587},
\seqnum{A104980},
\seqnum{A104980},
\seqnum{A111184},
\seqnum{A111529},
\seqnum{A111530},
\seqnum{A111531},
\seqnum{A111536},
\seqnum{A111544},
\seqnum{A111553}, and
\seqnum{A132159}.)

\end{document}